\documentclass[11pt]{amsart}
\usepackage[utf8]{inputenc}

\usepackage[english]{babel}
\usepackage[T1]{fontenc}

\usepackage{amsmath}
\usepackage{enumitem}
\usepackage{amssymb,amsfonts,stmaryrd,mathrsfs}
\usepackage{amsmath,latexsym,amsthm,mathtools,bbm}
\usepackage{pgf,tikz} 
\usepackage{subcaption}
\usepackage{hyperref}
\usepackage{todonotes}

\usepackage[capitalize]{cleveref}
\newtheorem{thm}{Theorem}[section]
\newtheorem{cor}[thm]{Corollary}
\newtheorem{lem}[thm]{Lemma}
\newtheorem{prop}[thm]{Proposition}
\newtheorem{defi}[thm]{Definition}
\newtheorem{conj}{Conjecture}

\theoremstyle{remark}
\newtheorem{rmq}{Remark}
\newtheorem{exe}{Example}[section]

\DeclareMathOperator{\bfx}{\mathbf{x}}
\DeclareMathOperator{\bfy}{\mathbf{y}}
\DeclareMathOperator{\bfxy}{\mathbf{xy}}
\DeclareMathOperator{\wv}{\mathcal{V}_\circ}
\DeclareMathOperator{\ver}{\mathcal{V}}
\DeclareMathOperator{\bv}{\mathcal{V}_\bullet}
\DeclareMathOperator{\cc}{cc}
\DeclareMathOperator{\hook}{hook}

\DeclareMathOperator{\type}{face-type}

\DeclareMathOperator{\wvt}{\vartheta_\circ^\rho}
\DeclareMathOperator{\F}{\mathit{F}^{(\alpha)}_{\lambda,\vartheta}}
\DeclareMathOperator{\Frho}{\mathit{F}^{(\alpha)}_{\lambda,\vartheta^\rho}}
\DeclareMathOperator{\B}{\mathit{B}^{(\alpha)}_{\mathit{n},\rho}}
\DeclareMathOperator{\Bj}{\mathit{B}^{(\alpha)}_{\mathit{j},\rho}}
\DeclareMathOperator{\MON}{MON}

\DeclareMathOperator{\D}{\mathcal{D}}
\DeclareMathOperator{\ID}{\mathcal{ID}}
\DeclareMathOperator{\IDO}{\mathcal{ID^O}}

 \address{Université de Lorraine\\ CNRS\\ IECL\\ F-54000 Nancy\\ France} 
\address{Université de Paris\\ CNRS\\ IRIF\\ F-75006\\ Paris\\ France}

\title{A note on the map expansion of Jack polynomials}
\author{Houcine Ben Dali}

\email{houcine.ben-dali@univ-lorraine.fr}

\begin{document}
\begin{abstract}
 
   In a recent work, Maciej Do\l{}e\k{}ga and the author have given a formula of the expansion of the Jack polynomial $J^{(\alpha)}_\lambda$ in the power-sum basis as a non-orientability generating series of bipartite maps whose edges are decorated with the boxes of the partition $\lambda$.  

    We conjecture here a variant of this expansion in which we restrict the sum on maps whose edges are injectively decorated by the boxes of $\lambda$. We prove this conjecture for Jack polynomials indexed by 2-column partitions. The proof uses a mix of combinatorial methods and differential operator computations.
\end{abstract}

\maketitle

\section{Introduction}

\subsection{Jack polynomials and maps}

Jack polynomials $J_\lambda^{(\alpha)}$ are symmetric functions indexed by an integer partition $\lambda$ and  a deformation parameter $\alpha$, and which have been introduced by Jack in \cite{Jack1970/1971}. Jack polynomials interpolate, up to scaling factors, between Schur functions for $\alpha=1$ and zonal polynomials for $\alpha=2$. 
In his work \cite{Stanley1989}, Stanley initiated the combinatorial analysis of these symmetric functions. Later, they have been connected to various objects of algebraic combinatorics, such as partitions, tableaux, paths and maps \cite{Macdonald1995,GouldenJackson1996a,KnopSahi1997,DolegaFeray2016,Moll2023,ChapuyDolega2022, BenDaliDolega2023}. 

%
In particular, Haglund and Wilson have given a combinatorial interpretation of the expansion of Jack polynomials in the power-sum basis in terms of weighted tableaux \cite{HaglundWilson2020}.  A different formula of this expansion has been obtained recently in \cite{BenDaliDolega2023} as a non-orientability generating series of maps whose edges are decorated with the boxes of a Young diagram.
The formula of \cite{BenDaliDolega2023} answers the long-standing question, going back to
to \cite{Hanlon1988,DolegaFeraySniady2014}, of giving an expression of Jack polynomials
in terms of maps, and it was used in \cite{BenDaliDolega2023} to prove positivity
conjectures on Jack characters. The purpose of this note is to study
variants of these formula, in which an additional injectivity property
is required. We make a conjecture in this direction, see \cref{conj injective}.

This injective conjecture is an $\alpha$-deformation of two known formula for Schur and Zonal functions (which correspond to Jack polynomials for $\alpha=1$ and $\alpha=2$ respectively). We prove \cref{conj injective} for Jack polynomials indexed by two-column partitions.


Roughly, a map is a graph drawn on a locally orientable surface. The study of maps is a well developed area with strong connections with analytic combinatorics, mathematical physics and probability \cite{BenderCanfield1986,LandoZvonkin2004,Chapuy2011,Eynard2016}. The relationship between generating series of maps and the theory of symmetric functions was first noticed via a character theoretic approach  \cite{JacksonVisentin1990,GouldenJackson1996a} and has then been developed to include other techniques such as matrix integrals and differential equations \cite{LaCroix2009,DolegaFeraySniady2014,ChapuyDolega2022,  BenDaliDolega2023}.

The introduction is organized as follows. In Sections \ref{ssec maps} and \ref{ssec decorated maps} we introduce some definitions related to maps.
In Section \ref{ssec non-inj}, we recall the maps expansion of Jack polynomials obtained in \cite{BenDaliDolega2023}. We conjecture an injective version of this expansion in \cref{ssec injective conj}. We formulate the main results of the paper in Sections \ref{ssec main thm} and \ref{ssec low degree}. In \cref{ssec related results}, we explain how our conjecture is connected to other conjectures relating Jack polynomials to maps. Finally, the outline of the paper is detailed in \Cref{ssec outline of the paper}.

\subsection{Maps}\label{ssec maps}
A \textit{connected map} is a connected graph embedded into a surface  such that all the connected components of the complement of the graph are simply connected (see \cite[Definition 1.3.7]{LandoZvonkin2004}). These connected components are called the \textit{faces} of the map. We consider maps up to homeomorphisms of the surface. A connected map is \textit{orientable} if the underlying surface is orientable. In this paper\footnote{This is not the standard definition of a map; usually a map is connected. }, a \textit{map} is an unordered collection of connected maps. A map is orientable if each one of its connected components is orientable. Finally, the \textit{size} of a map is its number of edges.

All maps considered here are \textit{bipartite}; \textit{i.e}  their vertices are colored in two colors, white and black, such that each edge connects two vertices of different colors. 
Note that in a bipartite map, all faces have even degree. 
Hence, we define the \textit{face-type} of a bipartite map of size $n$, as the partition of $n$ obtained by reordering the half degrees of the faces.
For a given map $M$, we denote by $|\ver(M)|$ its number of vertices, and by $\cc(M)$ its number of connected components. We also denote its number of white and black vertices by $|\wv(M)|$ and $|\bv(M)|$ respectively.

In order to enumerate maps with trivial automorphism group, we consider rooted maps; we say that a connected map is \textit{rooted} if it has  a distinguished white oriented corner $c$. The corner $c$ will be called \textit{the root corner} and  the edge following the root corner is called \textit{the root edge}.
More generally, a map is rooted if each one of its connected components is rooted. We say that a rooted map of size $n$ is \textit{labelled}, if its edges are labelled with the integers $1,2,...n$, such that the root edge of each connected component has the smallest label in its connected component. 

\subsection{\texorpdfstring{$\lambda$}{}-decorated maps}\label{ssec decorated maps}

\begin{defi}\label{def decorated maps}
We call a mapping $f$ of  $M$ in a Young diagram $\lambda$, a function which associates to each edge $e$ of $M$ a box of $\lambda$, such that
\begin{enumerate}
    \item if $e_1$ and $e_2$ are two edges incident to the same black vertex, then $f(e_1)$ and $f(e_2)$ are in the same row of $\lambda$.
    \item if $e_1$ and $e_2$ are two edges incident to the same white vertex, then $f(e_1)$ and $f(e_2)$ are in the same column of $\lambda$.
\end{enumerate}
Such pair $(M,f)$ is called a $\lambda$-decorated map. Moreover, we say that $(M,f)$ is a $\lambda$-injectively decorated maps, if $f$ satisfies the additional injectivity property:
\begin{enumerate} 
\setcounter{enumi}{2}
    \item two different edges are decorated by two different boxes; $f(e_1)\neq f(e_2)$ if $e_1\neq e_2$.
   
\end{enumerate}

We consider a total order on the boxes of $\lambda$ using the lexicographic order, see \cref{eq order} for a precise definition. A rooted $\lambda$-decorated map is a rooted map $M$ equipped with a function $f$ from $M$ to $\lambda$ with the condition that the root edge of each connected component is decorated with the smallest decorating box in its connected component.

We denote by $\D_n(\lambda)$ the set of $\lambda$-decorated maps of size $n$, orientable or not.
Similarly, we denote $\ID_n(\lambda)$ (resp. $\IDO_n(\lambda)$)  the set of rooted $\lambda$-injectively decorated maps of size $n$, orientable or not (resp. orientable).
\end{defi}

These definitions were introduced in \cite{FeraySniady2011a,FeraySniady2011} as injective mappings of permutations and matchings respectively into partition diagrams.    Later a reformulation in terms of maps has been considered in \cite{DolegaFeraySniady2014}. This reformulation is based on the fact that maps can be encoded using permutations or matchings (see \textit{e.g.} \cite{GouldenJackson1996a,LandoZvonkin2004}).

\begin{figure}[t]
\centering
\begin{subfigure}{.35\textwidth}
  \centering
    \includegraphics[width=0.3\textwidth]{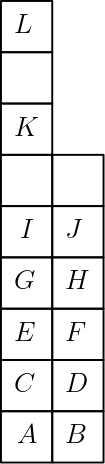}
    \subcaption{}
  \label{fig:sub1}
\end{subfigure}%
\begin{subfigure}{.65\textwidth}
  \centering
    \includegraphics[width=0.85 \textwidth]{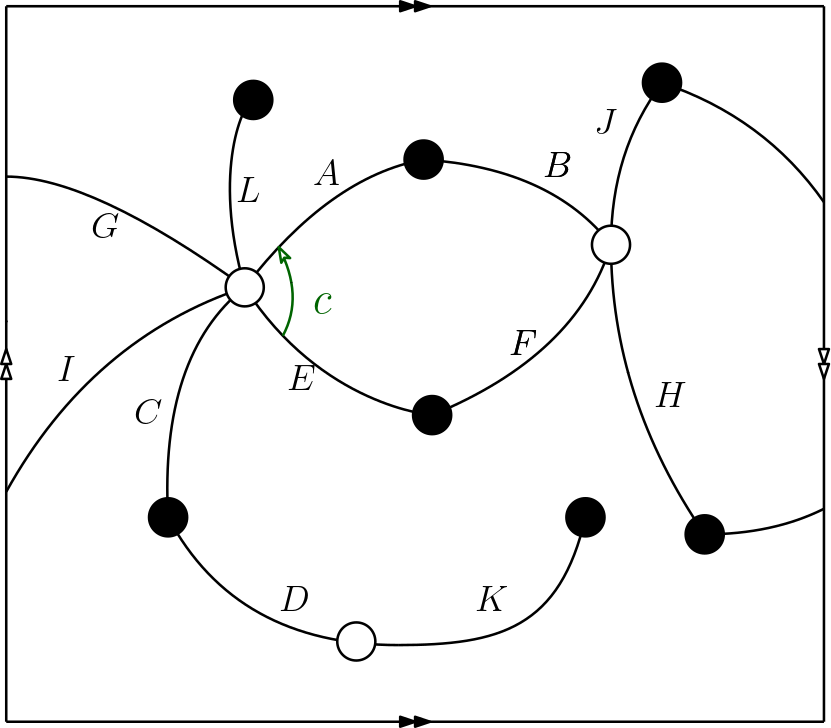}
    \subcaption{}
  \label{fig:sub2}
\end{subfigure}
\caption{The diagram of the partition $\lambda=2^61^3$ and a map in $\ID_{12}(\lambda)$ drawn the projective plan. The Latin letters indicate the decoration of the map edges with the boxes of the partition.}
\label{Figure map}
\end{figure}

\begin{exe}
We consider the 2-column partition $\lambda=2^61^3$ of size 15. In \cref{fig:sub1} we give an example of a rooted $\lambda$-injectively decorated map of size 12. In \cref{fig:sub2}, the left-hand side of the square should be glued to the right-hand
one (with a twist)  as indicated by the white arrows, and the top side should be glued to the
bottom one (without a twist) as indicated by the black arrows. The root corner is indicated by the green arrow.
\end{exe}

\subsection{The non-injective expansion of \texorpdfstring{\cite{BenDaliDolega2023}}{}}\label{ssec non-inj}

We start by the following definition due to Goulden and Jackson.


\begin{defi}[\cite{GouldenJackson1996}]
\textit{A statistic of non-orientability} on bipartite maps is a statistic $\vartheta$ with non-negative integer values, such that $\vartheta(M)=0$ if and only if $M$ is orientable.
\end{defi}

In practice, a statistic of non-orientability is supposed to "measure" the non-orientability of a map. Generally, such a measure is obtained by counting the number of "twisted" edges, following a given algorithm of decomposition of the map. Several examples of such statistics have been introduced in previous works \cite{LaCroix2009,DolegaFeraySniady2014,Dolega2017a,ChapuyDolega2022,BenDaliDolega2023}. 

In the following, maps are counted with a \textit{weight of non-orientability} of the form $b^{\vartheta(M)}$, where $b:=\alpha-1$ is the shifted Jack parameter. We denote by  $J^{(\alpha)}_\lambda(\mathbf{p})$ Jack polynomial associated to the partition $\lambda$ expressed in the power-sum basis $\mathbf{p}:=(p_1,p_2,...)$, see \cref{sec SymFun}.

\begin{thm}\cite{BenDaliDolega2023}\label{thm BDD23}
    Let $n$ be a positive integer and let $\lambda$ be a partition of $n$. There exists a statistic of non-orientability $\vartheta$ such that 
    $$  J_\lambda^{(\alpha)}(\mathbf{p})=\sum_{M\in\D_n(M)}(-1)^{n-|\wv(M)|}\frac{\alpha^{|\wv(M)|-\cc(M)}}{C(M)}b^{\vartheta(M)}p_{\type(M)} ,  
$$
where the sum is taken over maps in $\D_n(\lambda)$ considered with some labelling\footnote{The labelling used in \cite{BenDaliDolega2023} is different from the one considered here, and is related to the decoration of the map.}, and $C(M)$ is a normalization factor related to this labelling. 
\end{thm}

The first special cases of \cref{thm BDD23} which have been proved correspond to the cases $\alpha=1$ (Schur functions) and $\alpha=2$ (Zonal functions) see \cite{FeraySniady2011a, FeraySniady2011}. In both these papers, the authors start by proving the expansion in terms of maps injectively decorated. For $\alpha=1$ this corresponds to a well known formula for Schur functions \cite[Eq. (1.1)]{Hanlon1988}, (see also \cite[Proposition 5]{FeraySniady2011a} for a full proof of this result). We reformulate it here in the language of maps.
\begin{thm}\cite{Hanlon1988, FeraySniady2011a}\label{thm alpha=1}
    Let $\lambda$ be a partition of size $n$. Then 
\begin{equation}\label{eq alpha=1}
  J^{(\alpha=1)}_\lambda= \sum_{M\in\IDO_n(\lambda)}(-1)^{n-|\wv(M)|}p_{\type(M)}.  
\end{equation}
\end{thm}

For $\alpha=2$, a similar  injective formula has been proved in 
 \cite[Eq. (16)]{FeraySniady2011}. 

\begin{thm}\cite{FeraySniady2011}\label{Thm FS}
Let $\lambda$ be a partition of size $n$. Then 
\begin{equation}\label{eq alpha=2}
  J^{(\alpha=2)}_\lambda= \sum_{M\in\ID_n(\lambda)}(-1)^{n-|\wv(M)|}2^{|\wv(M)|-\cc(M)}p_{\type(M)}.  
\end{equation}
\end{thm}

Féray and \'{S}niady have proved in \cite{FeraySniady2011a} and \cite{FeraySniady2011}  that in the cases $\alpha=1$ and $\alpha=2$,  the injective and the non injective variants of these formulas are equivalent. To do so, they constructed a fixed points free involution on the set of maps equipped with a non injective mapping in the Young diagram of $\lambda$.  However, it seems hard to construct such an involution for general $\alpha$ because of the non-orientability weight. In \cite{BenDaliDolega2023}, the authors use an independent approach based on differential operators which allow to construct the generating series of decorated maps. They then prove that this series satisfies some characterization properties of the power-sum expansion of Jack polynomials, which allows to obtain directly the non-injective formula of \cref{thm BDD23}. Nevertheless, it seems natural to try to obtain an $\alpha$ deformation of the injective formulas of Equations \eqref{eq alpha=1} and \eqref{eq alpha=2}.

\subsection{An injective conjecture}\label{ssec injective conj}

We consider the following injective variant of \cref{thm BDD23}.
\begin{conj}\label{conj injective} 
Let $\lambda$ be a partition of size $n$.
Then there exists a statistic of non-orientability $\vartheta$ on $\ID_n(\lambda)$ such that, 
\begin{equation}\label{eq conj injective}
  J_\lambda^{(\alpha)}(\mathbf{p})=\sum_{M\in\ID_n(\lambda)}(-1)^{n-|\wv(M)|}\alpha^{|\wv(M)|-\cc(M)}b^{\vartheta(M)}p_{\type(M)} ,  
\end{equation}
\end{conj}

In the following, we will refer to the quantity  $(-1)^{|M|-|\wv(M)|}\alpha^{|\wv(M)|-\cc(M)}$ as the \textit{$\alpha$-weight} of the map $M$ and to $b^{\vartheta(M)}$ as its \textit{$b$-weight}.





 While the non injective formula of \cref{thm BDD23} allows to obtain polynomiality and positivity properties about Jack polynomials (see \cite[Theorem 1.6]{BenDaliDolega2023}), the injective conjecture we are considering here has the advantage of having less cancellation. Indeed, the fact that the considered maps are injectively decorated puts some restrictions on the underlying graph. For example, maps with multiple edges do not appear in \cref{eq conj injective}. This also implies that there exist no $\lambda$-injectively decorated maps of size bigger than $|\lambda|$ which would make obvious the vanishing property of Jack characters (see e.g \cite[Theorem 2.5]{BenDaliDolega2023}).

Both the particular cases of  \eqref{eq alpha=1} and \eqref{eq alpha=2} are obtained using representation theory tools which seem to be specific to these specialisations and can not been used for general $\alpha$. New tools are then required to understand \cref{conj injective}.

\subsection{Main theorem}\label{ssec main thm}
The main contribution of this paper is to establish \cref{conj injective} for Jack polynomials indexed by 2-column partitions.
\begin{thm}\label{main thm}
Let $\lambda$ be a 2-column partition of size $n$.
There exists an explicit statistic of non-orientability $\vartheta$ on $\ID_n(\lambda)$ such that 
\begin{equation}\label{eq main thm}
J_\lambda^{(\alpha)}(\mathbf{p})=\sum_{M\in\ID_n(\lambda)}(-1)^{n-|\wv(M)|}\alpha^{|\wv(M)|-\cc(M)}b^{\vartheta(M)}p_{\type(M)}.
\end{equation}
\end{thm}

 Unfortunately, the statistic $\vartheta$ that we use here for 2-column partitions does not work for all partitions $\lambda$.

We now briefly describe  our proof strategy.
The starting point is an expression of the  Jack polynomial $J_\lambda^{(\alpha)}$, indexed by a 2-column partition $\lambda=[2^r,1^s]$ in the monomial basis which is due to Stanley \cite[Equation 46]{Stanley1989}:
\begin{equation}\label{eq 2-column Jack}
J_\lambda^{(\alpha)}(\mathbf{p})=\sum_{i=1}^r(r)_i(\alpha+r+s)_i(2(r-i)+s)!m_{2^i1^{2(r-i)+s}}(\bfx),
\end{equation}
where we denote the falling factorial by
$$(a)_k:=a(a-1)...(a-k+1),$$
for any real $a$, and any non negative integer $k$.

On the other hand, using a variant of the result of Chapuy and Do\l{}e\k{}ga \cite{ChapuyDolega2022}, we give a particularly simple expression in the monomial basis for a generating series of labelled bipartite maps with non-orintability weights, which are not yet decorated, see \cref{thm bipartite maps}.
We prove then that the generating series of $\lambda$-injectively decorated maps can be obtained from the series of labelled maps using an "leaf addition" procedure, see \cref{prop F} and \cref{lem X+ Y+}.

After some computations, we show that the right-hand side of \cref{eq main thm} is equal to that of \cref{eq 2-column Jack}, concluding the proof.

\subsection{Low degree terms in the parameter \texorpdfstring{$\alpha$}{}}\label{ssec low degree}
It is known from the work of Lapointe and Vinet \cite{LapointeVinet1995}, that the coefficients of the Jack polynomials in the power-sum basis are polynomial in the parameter $\alpha$. 
As a consequence of the main result, we prove in \cref{sec applications} that \cref{conj injective} holds when we extract the coefficients of $\alpha^0$ and $\alpha^1$ in \cref{eq conj injective} for any partition $\lambda$.

\subsection{Connection to other conjectures relating maps to Jack polynomials}\label{ssec related results}

Hanlon has conjectured in \cite{Hanlon1988} that Jack polynomials expansion in the power-sum basis can be expressed as  weighted sums of pairs of permutation. This conjecture can be reformulated as follows.
\begin{conj}\cite{Hanlon1988}\label{conj hanlon}
Fix a partition $\lambda$ of size $n$. 
Then, there exists a statistic $w$ on orientable $\IDO_n(\lambda)$ maps such that 
\begin{equation}\label{eq inj conj}
J^{(\alpha)}_\lambda(\mathbf{p})=\sum_{M\in\IDO_n(\lambda)}(-1)^{n-|\wv(M)|}\alpha^{w(M)}p_{\type(M)}.
\end{equation}
\end{conj}
Hence, \cref{conj injective} can be seen as a non orientable version of Hanlon's conjecture. 
Indeed, considering Theorems \ref{Thm FS} and \ref{thm BDD23} above and other works relating Jack polynomials to maps
\cite{GouldenJackson1996,ChapuyDolega2022}, it seems more natural to take a sum over non-oriented maps (as we do in \cref{conj injective}) rather than a sum over oriented maps (as was done by Hanlon in \cref{conj hanlon}).


Note that for $\alpha=1$ both conjectures are equivalent to \cref{thm alpha=1}.
An other special case for which the two conjectures are equivalent, is the case of hook-shaped partitions $\lambda$, \textit{i.e} $\lambda=[r,1^s]$ for some integers $r$ and $s$. Indeed, it follows from the definition of decorated maps that only orientable maps appear in the sum of \cref{conj injective}.
This case has been treated  by Hanlon \cite[Property 2]{Hanlon1988}, see also \cite[Example 4.(b) page 385]{Macdonald1995}.



Moreover, both the result of \cite{BenDaliDolega2023} and \cref{conj injective} can be seen as dual questions for the celebrated Matchings-Jack and $b$-conjectures of Goulden and Jackson, which suggest that the multivariate non-orientability generating series of bipartite maps can be expressed through Jack polynomials. These conjectures are still open despite many partial results \cite{DolegaFeray2016,DolegaFeray2017,Dolega2017a,KanunnikovPromyslovVassilieva2018,ChapuyDolega2022,BenDali2022a,BenDali2023}.

\subsection{Outline of the paper}\label{ssec outline of the paper}
In \cref{sec preliminaries}, we introduce some useful notation and preliminary results. Based on the work of Chapuy and Do\l{}e\k{}ga, we give in \cref{sec generating series} an expression of a particular specialization of the generating series of weighted labelled maps in terms of monomial functions. \cref{sec proof main thm} is dedicated to the proof of \cref{main thm}. As a consequence of the main result, we establish in \cref{sec applications} the \cref{conj injective} for terms of degree 0 and 1 in $\alpha$.

\section{Notation and preliminaries}\label{sec preliminaries}
For the definitions and notation introduced in \cref{subsec Partitions,sec SymFun} 
we refer to \cite{Stanley1989,Macdonald1995}.
\subsection{Partitions}\label{subsec Partitions}

A \textit{partition} $\lambda=[\lambda_1,...,\lambda_\ell]$ is a weakly decreasing sequence of positive integers $\lambda_1\geq...\geq\lambda_\ell>0$. The integer $\ell$ is called the \textit{length} of $\lambda$ and is denoted $\ell(\lambda)$. The size of $\lambda$ is the integer $|\lambda|:=\lambda_1+\lambda_2+...+\lambda_\ell.$ If $n$ is the \textit{size} of $\lambda$, we say that $\lambda$ is a partition of $n$ and we write $\lambda\vdash n$. The integers $\lambda_1$,...,$\lambda_\ell$ are called the \textit{parts} of $\lambda$. For every $i\geq1$, we denote by $m_i(\lambda)$ the number of parts equal to $i$ in $\lambda$, and we introduce the following notation
$$z_\lambda:=\prod_{i\geq1}m_i(\lambda)!i^{m_i(\lambda)}.$$ 
We recall the \textit{dominance partial} ordering on partitions $\leq$ defined by 
$$\mu\leq\lambda \iff |\mu|=|\lambda| \text{ and }\hspace{0.3cm} \mu_1+...+\mu_i\leq \lambda_1+...+\lambda_i \text{ for } i\geq1.$$



\noindent We identify a partition  $\lambda$ with its Young diagram defined by 
$$\lambda:=\{(i,j),1\leq i\leq \ell(\lambda),1\leq j\leq \lambda_i\}.$$

For a given partition $\lambda$, we define the total order on the boxes of $\lambda$ given by the lexicographic order on their coordinates: 
\begin{equation}\label{eq order}
    \Box_1=(i_1,j_1)< \Box_2=(i_2,j_2) \Longleftrightarrow i_1<i_2 \text{ or } (i_1=i_2 \text{ and }j_1<j_2).
\end{equation}
\textit{The conjugate partition} of $\lambda$, denoted $\lambda^t$, is the partition associated to the Young diagram obtained by reflecting the diagram of $\lambda$ with respect to the line $j=i$:
$$\lambda^t:=\{(i,j),1\leq j\leq \ell(\lambda),1\leq i\leq \lambda_i\}.$$

Fix a box $\Box:=(i,j)\in\lambda$. Its \textit{arm-length} is given by 
$$a_\lambda(\Box):=|\{(i,r)\in\lambda,r>j\}|=\lambda_i-j,$$ 
and its \textit{leg-length} is given by 
$$\ell_\lambda(\Box):=|\{(r,j)\in\lambda,r>i\}|=(\lambda^t)_j-i.$$ 

Two $\alpha$-deformations of the hook-length product were introduced in \cite{Stanley1989};
$$\hook_\lambda^{(\alpha)}:=\prod_{\Box\in\lambda}\left(\alpha a_\lambda(\Box)+\ell_\lambda(\Box)+1\right),
\qquad 
\hook_\lambda'^{(\alpha)}:=\prod_{\Box\in\lambda}\left(\alpha(a_\lambda(\Box)+1)+\ell_\lambda(\Box)\right).$$
With these notation, the classical hook-length product is given by 
$$H_\lambda:=\hook_\lambda^{(1)}=\hook_\lambda'^{(1)}.$$
Finally, we define the \textit{$\alpha$-content} of a box $\Box:=(i,j)$ by $c_\alpha(\Box):=\alpha(j-1)-(i-1)$.

\subsection{Symmetric functions and Jack polynomials}\label{sec SymFun}

We fix an alphabet $\mathbf{x}:=(x_1,x_2,..)$. We denote by $\mathcal{S}$ the algebra of symmetric functions in $\mathbf{x}$ with coefficients in $\mathbb Q$. For every partition $\lambda$, we denote $m_\lambda$ the monomial function and $p_\lambda$ the power-sum function associated to the partition $\lambda$.
We  consider the associated alphabet of power-sum functions $\mathbf{p}:=(p_1,p_2,..)$. 

Let $\mathcal{S}_\alpha:=\mathbb{Q}[\alpha]\otimes\mathcal{S}$ the algebra of symmetric functions with rational coefficients in $\alpha$. 
We denote by $\langle.,.\rangle_\alpha$ the $\alpha$-deformation of the Hall scalar product defined by 
$$\langle p_\lambda,p_\mu\rangle_\alpha=z_\lambda\alpha^{\ell(\lambda)}\delta_{\lambda,\mu},\text{ for any partitions }\lambda,\mu.$$

Macdonald \cite[Chapter VI.10]{Macdonald1995} has proved that there exists a unique family of symmetric functions $(J_\lambda^{(\alpha)})$ in $\mathcal{S}_\alpha$ indexed by partitions, satisfying the following properties:
$$\left\{ \begin{array}{ll}
       \text{Orthogonality: }
     &\langle J_\lambda,J_\mu\rangle_\alpha=0, \text{ for }\lambda\neq\mu,\\
      \text{Triangularity: }
     &[m_\mu]J_\lambda=0, \text{ unless }\mu\leq\lambda,\\
     \text{Normalization: }
     &[m_{1^n}]J_\lambda=n!, \text{ for }\lambda\vdash n,
\end{array}\right.$$
where $[m_\mu]J_\lambda$ denotes the coefficient of $m_\mu$ in $J_\lambda$, and $1^n$ is the partition with $n$ parts equal to~1. These functions are known as the \textit{Jack polynomials}.
In particular, Jack polynomials indexed by 1-column partitions are given by
\begin{equation}\label{eq column Jack}
    J^{(\alpha)}_{1^n}=n!m_{1^n}=\sum_{\mu\vdash n}(-1)^{n-\ell(\mu)}\frac{n!}{z_\lambda}p_\mu,
\end{equation}
and are independent of the parameter $\alpha$.
The squared norm of Jack polynomials can be expressed in terms of the deformed hook-length products, see \cite[Theorem 5.8]{Stanley1989}:
\begin{equation}\label{eq j alpha}
    j_\lambda^{(\alpha)}:=\langle J_\lambda,J_\lambda\rangle_\alpha=\hook^{(\alpha)}_\lambda\hook'^{(\alpha)}_\lambda.
\end{equation}
In particular, we have
\begin{equation}\label{eq column j}
j^{(\alpha)}_{1^n}=n!(\alpha+n-1)_n.    
\end{equation}
In this paper, Jack polynomials will always be expressed in the power-sum variables $\mathbf{p}$ rather than the alphabet $\bfx$ (this is possible since the power-sum functions form a basis of the symmetric functions algebra).

We have the following theorem due to Macdonald \cite[Chapter VI Eq. 10.25]{Macdonald1995}, which gives an expression of Jack polynomials when all power-sum variables $p_i$ are equal to a variable $u$.
\begin{thm}[\cite{Macdonald1995}]\label{Jack formula}
For every $\lambda\in\mathcal{P}$, we have
$$J_\lambda^{(\alpha)}(\underline{u})=\prod_{\Box\in\lambda}\left(u+c_\alpha(\Box)\right),$$
where $\underline{u}:=(u,u,...)$. 
\end{thm}
We conclude this subsection with the following corollary. 
\begin{cor}\label{cor ev -alpha}
Let $\lambda\vdash n\geq1.$
We have the following expression for Jack polynomials specialized at $u=-\alpha$:
$$J_\lambda^{(\alpha)}(\underline{-\alpha})=\left\{
\begin{array}{cc}
    (-1)^n (\alpha+n-1)_n  & \text{if } \lambda=1^n, \\
    0 & \text{otherwise.}
\end{array}
\right.$$
\end{cor}

\subsection{Side-marked maps}\label{ssec side-marked}
Recall that a connected map is rooted if it has a marked oriented white corner. This is equivalent to saying that it has a marked edge-side (by convention this will be the side following the root corner).
We say that a map is \textit{side-marked} if each one of its edges has a distinguished side.
Rather than using rooted maps, in which only one edge in each connected component has a distinguished side, it will be more practical in the proof of the main result to consider side-marked maps.
Note that to a rooted map of size $n$, we can associate $2^{n-\cc(M)}$ different side-marked maps by choosing the distinguished sides of non-root edges. 

In this paper, we consider statistics of non-orientability on labelled and injectively decorated maps. Such statistics only depend on the labelling or the decoration of the map, but not on its  rooting or side-marking. 

\section{Generating series of \texorpdfstring{$b$}{}-weighted bipartite maps}\label{sec generating series}

In this section, we introduce a family of non-orientability statistics  $\wvt$ on bipartite maps.
We consider then the generating series of labelled bipartite maps, where each map $M$ is counted with a $b$-weight $b^{\wvt(M)}$, and we use a result of Chapuy and Do\l{}e\k{}ga to give an expression of this function in terms of monomial functions. 

\subsection{A statistic of non-orientability on labelled map}\label{ssec stat}
\begin{figure}[t]
\centering

\begin{subfigure}{0.8\textwidth}
\centering

\begin{subfigure}{.4\textwidth}
  \centering
    \includegraphics[width=0.5\textwidth]{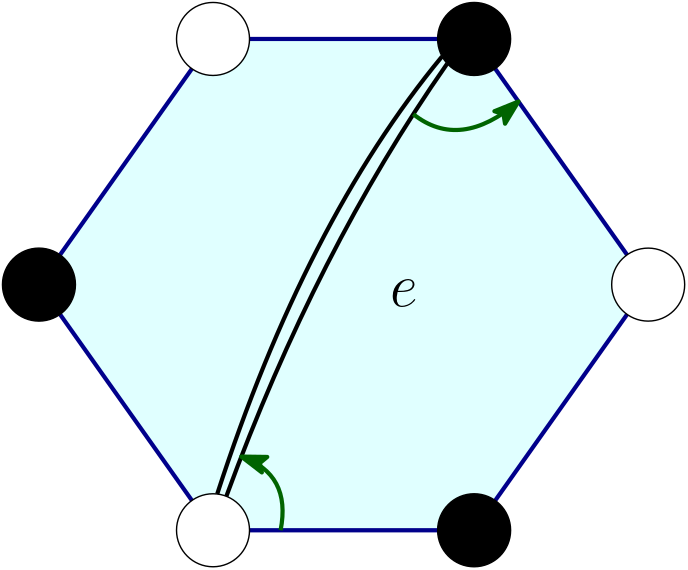}
    \label{diagonal}
\end{subfigure}%
\begin{subfigure}{.4\textwidth}
  \centering
    \includegraphics[width=0.5 \textwidth]{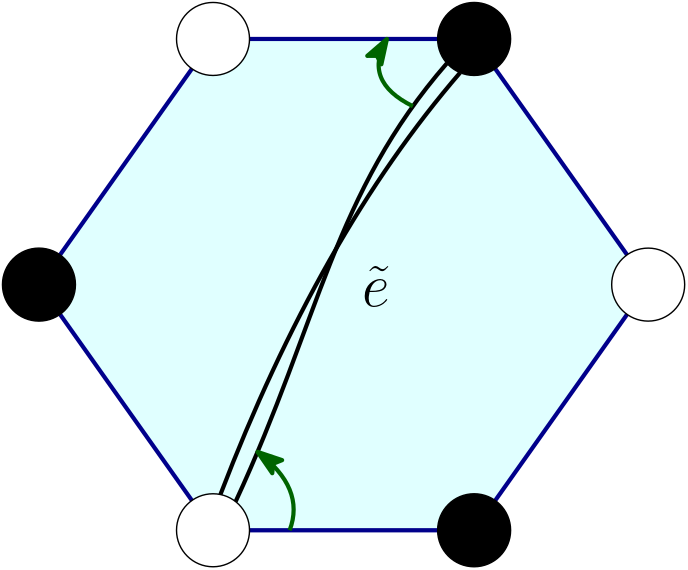}
\label{twist}
\end{subfigure}

\label{diagonal twist}
\caption{A pair of twisted edges $(e,\tilde{e})$ between two corners of the same face; $e$ is a diagonal while $\tilde{e}$ is a twist. }
\end{subfigure}

\begin{subfigure}{0.8\textwidth}
\begin{subfigure}{.45\textwidth}
  \centering
    \includegraphics[width=0.9\textwidth]{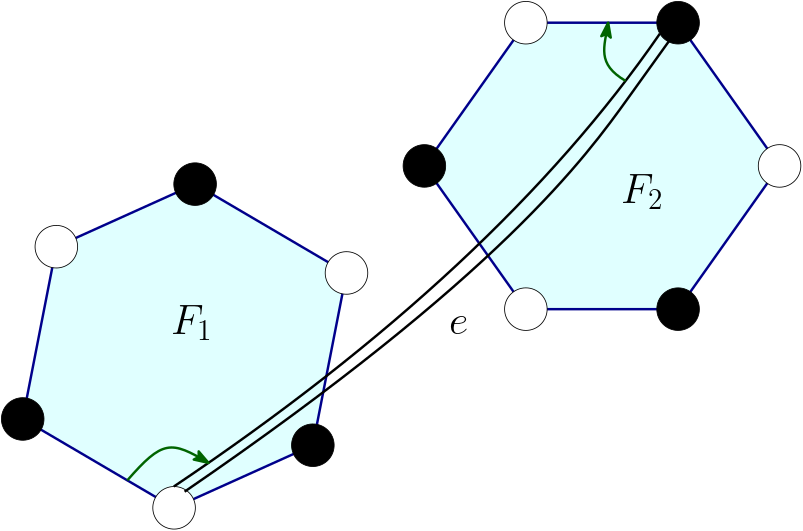}
  \label{handle 1}
\end{subfigure}%
\begin{subfigure}{.45\textwidth}
  \centering
    \includegraphics[width=0.9 \textwidth]{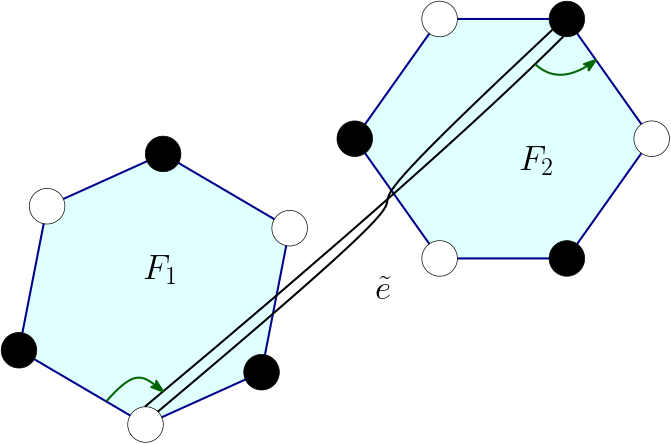}
    \label{handle 2}
\end{subfigure}
 \caption{A pair of twisted edges $(e,\tilde{e})$ between two corners of different faces $F_1$ and $F_2$. If $F_1$ and $F_2$ are on the same connected component then the edges $e$ and $\tilde e$ are handles, otherwise they are bridges. }
\label{handles bridges}
 \end{subfigure}

\caption{The different ways of adding an edge on a map. The added edge is represented each time by a band, that can be twisted at most once. The arrows indicate how the added edge connects the two respective corners.}
\label{fig edge-types}
\end{figure}

Let $M$ be a bipartite map and let $c_1$ and $c_2$ be two corners of $M$ of different colors. Then we have two ways to add an edge to $M$ between these two corners. We denote by $e_1$ and $e_2$ these edges. We say that the pair $(e_1,e_2)$ is a \textit{pair of twisted edges} on the map $M$ and we say that $e_2$ is obtained by twisting $e_1$ (see \cref{fig edge-types}). Note that if $M$ is connected and orientable, then exactly one of the maps $M\cup\{e_1\}$ and $M\cup\{e_2\}$ is orientable. 
For a given map with a distinguished edge $(M,e)$, we denote $(\widetilde{M},\tilde{e})$ the map obtained by twisting the edge $e$.  

We recall that $b$ is the parameter related to the Jack parameter $\alpha$ by $b=\alpha-1$.
We now give the definition of a measure of non-orientability due to La Croix.

\begin{defi}\cite[Definition 4.1]{LaCroix2009}
We call a measure of non-orientability \textup{(MON)} a function $\rho$ defined on the set of maps $(M,e)$ with a distinguished edge, with values in  $\left\{1,b\right\}$, satisfying the following conditions:
\begin{itemize}
    \item if $e$ connects two corners of the same face of $M\backslash\{e\}$, and the number of the faces increases by 1 by adding the edge $e$ on the map $M\backslash\{e\}$, then $\rho(M,e)=1$. In this case we say that $e$ is a \textit{diagonal}.
    \item  if $e$ connects two corners in the same face $M\backslash\{e\}$, and the number of the faces of $M\backslash\{e\}$ is equal to the number of faces of $M$, then $\rho(M,e)=b$. In this case we say that $e$ is a \textit{twist}.
    \item if $e$ connects two corners of two different faces lying in the same connected component of $M\backslash\{e\}$, then $\rho$  satisfies $\rho(M,e)+\rho(\widetilde{M},\tilde{e})=1+b$. Moreover, if $M$ is orientable then $\rho(M,e)=1$. In this case we say that $e$ is a \textit{handle}.
    \item if $e$  connects two faces lying in two different connected components, then $\rho(M,e)=1$. In this case, we say that $e$ is a \textit{bridge}.
\end{itemize}

\end{defi}

We now introduce the following statistic of non-orientability on labelled maps which will have a key role in this paper.
\begin{defi}\label{def maps decomposition}
Let $\rho$ be a \textup{MON} and let $M$ be a connected labelled map of size $n$. We define the $b$-weight $\rho(M)$ of $M$, as the weight obtained by decomposing $M$ by deleting the edges in a decreasing order of their labels;
$$\rho(M):=\prod_{1\leq j\leq n}\rho(M_j,e_j),$$
where $e_j$ is the edge of $M$ labelled by $j$, and $M_j$ is the map obtained by deleting the edges $e_n$, $e_{n-1}$,..., $e_{j+1}$ from $M$. 
We extend  this definition by multiplicativity for disconnected maps;  we define the $b$-weight of a labelled disconnected map $M$ by 
$$\rho(M):=\prod_i\rho(M_i),$$
where the product is taken over the connected components of $M$. We define the statistic $\vartheta_\circ^\rho$ on labelled bipartite maps by $\rho(M)=b^{\wvt(M)}$. 
\end{defi}


It follows from the definitions that $\wvt$ is a statistic of non-orientability over labelled bipartite maps.
\begin{rmq}\label{rq statistics}
The statistic $\wvt$ is a variant of the statistics introduced in \cite{LaCroix2009,Dolega2017a,ChapuyDolega2022}; in these papers an "edge-deletion" procedure is used to decompose the maps while we use here the order given by the labels of the edges.
\end{rmq}

\subsection{Generating series of labelled bipartite maps}
Let $u_1$ and $u_2$ be two variables. 
We introduce the generating series of labelled rooted bipartite maps of size $n$:
\begin{equation}\label{eq B}
\B(\bfx,u_1,u_2)= \sum_{M}\frac{1}{n!}u_1^{|\wv(M)|}u_2^{|\bv(M)|}\alpha^{-\cc(M)}b^{\wvt(M)}p_{\type(M)}(\bfx),    
\end{equation}
where the sum is taken over labelled rooted bipartite maps of size $n$, and $p_\mu$ is the power-sum function associated to the partition $\mu$.

The following theorem is a variant of a special case in \cite[Theorem 5.10]{ChapuyDolega2022}, where a different statistic of non-orientability is used (see \cite[Definition 3.6]{ChapuyDolega2022} and \cref{rq statistics}). For completeness, we give more details about the proof of this theorem in \cref{appendix}.
\begin{thm}[\cite{ChapuyDolega2022}]\label{thm CD20}
For every $n\geq 1$, we have 
$$\B(\bfx,u_1,u_2)=\sum_{\xi\vdash n}\frac{J^{(\alpha)}_\xi(\mathbf{p})J^{(\alpha)}_\xi(\underline{u_1)}J^{(\alpha)}_\xi(\underline{u_2)}}{j^{(\alpha)}_\xi},$$
where $\mathbf{p}$ denotes the power-sum alphabet related to $\mathbf{x}$ \textup{(}see \cref{sec SymFun}\textup{)}.
\end{thm}

We deduce the following corollary which will be useful in the proof of \cref{main thm}.

\begin{cor}\cite{ChapuyDolega2022}\label{thm bipartite maps}
For every $n\geq0$ and for every \textup{MON} $\rho$, we have
$$B^{(\alpha)}_{n,\rho}(\bfx,-\alpha,-\alpha)=\frac{(\alpha+n-1)_n}{n!}J^{(\alpha)}_{1^n}(\mathbf{p})=(\alpha+n-1)_nm_{1^n}(\bfx).$$

\begin{proof}
Specializing $u_1=u_2=-\alpha$ in \cref{thm CD20} and using \cref{cor ev -alpha} we get that 
$$B^{(\alpha)}_{n,\rho}(\bfx,-\alpha,-\alpha)=\frac{\left((\alpha+n-1)_n\right)^2}{j^{(\alpha)}_\xi}J^{(\alpha)}_{1^n}(\mathbf{p}).$$
We conclude using \cref{eq column j,eq column Jack}.
\end{proof}
\end{cor}

\section{Proof of \texorpdfstring{\cref{main thm}}{}}\label{sec proof main thm}
\subsection{Preliminaries}

For every statistic of non-orientability $\vartheta$, we introduce  the generating series of $\lambda$-injectively decorated bipartite maps:

\begin{align*}
    \F(\mathbf{p)}:&=\sum_{\text{rooted maps }M}(-1)^{n-|\wv(M)|}\alpha^{|\wv(M)|-\cc(M)}b^{\vartheta(M)}p_{\type(M)}\\
    &=\sum_{\text{side-marked maps }M}\frac{(-1)^{n-|\wv(M)|}}{2^{|\lambda|-\cc(M)}}\alpha^{|\wv(M)|-\cc(M)}b^{\vartheta(M)}p_{\type(M)},
\end{align*}
where the two sums run over $\lambda$-injectively decorated  bipartite maps $M$ of size $|\lambda|$, which are rooted in the first line and side-marked in the second one, see \cref{ssec side-marked}. 
Hence \cref{main thm} can be reformulated as follows:
for some statistic of non-orientability $\vartheta$ we have $J_\lambda^{(\alpha)}=\F$.

The idea of the proof is to expand $\F$ in the monomial basis and compare the expression obtained to \cref{eq 2-column Jack}. To this purpose, we rewrite $\F$ as a sum on colored maps;
we call a \textit{coloring} of bipartite map a function $\mathcal{C}$ on the faces of $M$ with positive integer values. We say then that $(M,\mathcal{C})$ is \textit{a colored map}. For any colored map $(M,\mathcal{C})$, we define 
\begin{equation*}
x^{(M,\mathcal{C})}:=\prod_{f}x_{\mathcal{C}(f)}^{\deg(f)},    
\end{equation*}
where the product is taken over the faces $f$ of $M$.
Hence,
\begin{equation}\label{eq p x}
p_{\type(M)}(\mathbf{x})=\sum_{\mathcal{C}}x^{(M,\mathcal{C})},    
\end{equation}
where the sum is taken over all the colorings of $M$.
With this notation, the generating series $\F$ has the following expression in the alphabet $\bfx$:
$$\F(\bfx):=\sum_{(M,\mathcal{C})}\frac{(-1)^{n-|\wv(M)|}}{2^{|\lambda|-\cc(M)}} \alpha^{|\wv(M)|-\cc(M)}b^{\vartheta(M)}x^{(M,\mathcal{C})},$$
where the sum runs over colored $\lambda$-injectively decorated side-marked maps of size $|\lambda|$.

\subsection{From \texorpdfstring{$\lambda$}{}-decorated maps to labelled \texorpdfstring{$\circ^1/\circ^2$}{}-bipartite maps} \label{ssec decorated labelled maps}

\begin{defi}
Fix a  2-column  partition $\lambda$ and let $M$ be a $\lambda$-injectively decorated map. 
We say that a white vertex $v$ of $M$ has color $\circ^1$ \textup{(}resp. $\circ^2$\textup{)} if the edges incident to $v$ have labels in the first \textup{(}resp. the second\textup{)} column of $\lambda$.
For a face $f$ of $M$, we denote by $\deg_{\circ^1}(f)$ \textup{(}respectively $\deg_{\circ^2}(f)$\textup{)} the number of white corners of color $\circ^1$ \textup{(}respectively $\circ^2$\textup{)} incident to $f$.
\end{defi}

Let $\lambda$ be a 2-column partition, and let $M$ be a $\lambda$-injectively decorated map. Note that this implies that all black vertices of $M$ have degree 1 or 2. 
We call a \textit{black-leaf edge} of $M$ an edge incident to a black vertex of degree 1. 
We define the labelled bipartite map $M_\circ$ obtained by forgetting the black vertices $v$ of $M$ as follows:
\begin{itemize}
    \item If $v$ has degree 1, then we delete $v$ and the edge incident to it. Hence, we delete all the black-leaf edges of $M$. 
    \item If $v$ has degree 2, then we forget the vertex $v$ and consider the two edges incident to it as one edge (hence this edge separates two vertices of color $\circ^1$ and $\circ^2)$. 
\end{itemize}

Note that the map $M_\circ$ hence obtained has only white vertices, and is bipartite with respect to the colors $\circ^1$ and $\circ^2$. Moreover, we have an injection from the edges of  $M_\circ$ to the rows of size 2 of $\lambda$. Hence $M_\circ$ comes with a natural labelling inherited from the row indices associated to the edges via this injection. 

We now define a statistic of non-orientability on $\lambda$-injectively decorated maps.

\begin{defi}[A statistic on $\lambda$-injectively decorated maps]\label{def stat vartheta}
Fix a 2-column partition $\lambda$.
For every MON $\rho$, we associate to each $\lambda$-injectively decorated map $M$ the statistic $\vartheta^\rho(M):=\wvt(M_\circ)$, where $\wvt$ is the statistic on labelled bipartite maps defined in \cref{ssec stat}, and $M_\circ$ is the labelled bipartite map obtained from $M$ by forgetting the black vertices as explained above.
\end{defi}

Conversely, if $M_\circ$ is a labelled side-marked map of size $j\leq r$ which is bipartite in the colors $\circ^1$ and $\circ^2$, we obtain a $\lambda$-injectively decorated side-marked map of size $|\lambda|$ by realizing the following steps:
\begin{enumerate}
    \item  We choose a set $I$ of $j$ rows of $\lambda$ of size 2 (we have $\binom{r}{j}$ ways to choose such a set), and we associate to each edge of $M$ a row in $I$, using the labelling of $M$. 
    \item We add a black vertex in the middle of each edge.  Hence, we transform an edge $e$ to two edges $e_1$ and $e_2$, that connects respectively a black vertex to white vertices of color $\circ^1$ and $\circ^2$. Notice that we multiply then the size of the map by 2.
    \item  If $e$ is associated to a row $R_e$ of size $2$ in $\lambda$ then we decorate the edge $e_1$ (resp. $e_2$) by the box of $R_e$ which is in the first column (resp. the second column).
    \item For each box $\Box$ of $\lambda$ which is not used in the decoration of the map, we add a black-leaf edge connected to a white vertex of color $i$ (possibly a new white vertex), where $i\in\{1,2\}$ is the column containing $\Box$. Moreover, we add these edges successively in a increasing order of the decorating boxes. Note that each time we add a black-leaf edge connected to an existing white vertex then there are two different ways to mark one of its sides.
\end{enumerate}
 Note that each connected component of $M$ inherits a root from $M_\circ$ which is given by a an oriented corner of color $\circ^1$. 

\subsection{Adding black leaves}
\noindent We consider a second alphabet $\bfy:=(y_1,y_2,..)$ and we define the product alphabet $\bfxy:=(x_1y_1,x_2y_2,...)$. 
We introduce the two following operators 
$$X_+:=\sum_{i\geq1}\left(x_i-x_i^2\frac{\partial}{\partial x_i}\right),\qquad
Y_+:=\sum_{i\geq1}\left(y_i-y_i^2\frac{\partial}{\partial y_i}\right).$$
Let $(M,\mathcal{C})$ be a colored $\lambda$-injectively decorated side-marked bipartite map. We define the marking of  $(M,\mathcal{C})$ by:
$$\kappa(M,\mathcal{C}):=\frac{(-1)^{n-|\wv(M)|}}{2^{|\lambda|-\cc(M)}}\alpha^{|\wv(M)|-\cc(M)}b^{\vartheta^\rho(M)}\prod_{f}\left[x_{\mathcal{C}(f)}^{\deg_{\circ^1}(f)}y_{\mathcal{C}(f)}^{\deg_{\circ^2}(f)}\right],$$
where the product is taken over all the faces of $M$.  Hence the alphabet $\bfx$ (resp. $\bfy$) is a marking for $\circ^1$ (resp. $\circ^2$) colored corners.

Fix two $\lambda$-injectively decorated maps $M$ and $N$, such that $M$ is obtained from $N$ by deleting some black-leaf edges. 
If the map $N$ is equipped with a coloring $\mathcal{D}$, then $M$ is naturally equipped with a coloring $\mathcal{C}$, where we use the convention that deleting a black vertex of degree 1 from a face does not change its color. We say that the coloring $\mathcal{C}$ is \textit{inherited} from the coloring $\mathcal{D}$. 

The operator $X_+$ and $Y_+$ allow to add  black-leaf edges incident respectively to a white vertex of color $\circ^1$ and $\circ^2$. More precisely, we have the following lemma.
\begin{lem}\label{lem add black leaves}
Fix a $2$-column partition  $\lambda=2^r1^s$ and 
let $(M,\mathcal{C})$ be a colored $\lambda$-injectively decorated side-marked map.
Then,
\begin{equation}\label{eq X+}
     X_+\kappa(M,\mathcal{C})=\sum_{(M\cup{e},\widetilde{\mathcal{C}})}\kappa(M\cup e,\widetilde{\mathcal{C}}).
\end{equation}
\begin{equation}\label{eq Y+}
    Y_+\kappa(M,\mathcal{C})=\sum_{(M\cup{e},\widetilde{\mathcal{C}})}\kappa(M\cup e,\widetilde{\mathcal{C}}).
\end{equation}
where the sum in \cref{eq X+} \textup{(}resp. \cref{eq Y+}\textup{)} is taken over all $\lambda$-injectively decorated side-marked maps obtained by adding a  black-leaf edge $e$ to $M$ such that
\begin{itemize}
    \item the  colored map $(M\cup{e},\widetilde{\mathcal{C}})$ obtained is such that $\mathcal{C}$ is inherited from $\widetilde{\mathcal{C}}$.
    \item the edge $e$ is decorated by the smallest box in the first column \textup{(}resp. second column\textup{)}, which is not yet used in the decoration of $M$.
\end{itemize}
 
\begin{proof}
Let us prove \cref{eq X+}.
We start by noticing that adding a black-leaf edge to a map does not change its $b$-weight, this is straightforward from \cref{def stat vartheta}. 
We have two ways to add a black leaf-edge decorated by a box in the first column of $\lambda$:
\begin{itemize}
    \item by adding an isolated edge.
    \item by adding a black leaf on a $\circ^1$-corner.
\end{itemize}

In the first case the $\alpha$-weight of the map does not change, since the size of the map, the number of white vertices and the number of connected components all increase by 1. Finally, we choose a color $i\geq1$ for the new face, this is guaranteed by the operator $\sum_{i\geq1} x_i$. 

In the second case, the $\alpha$-weight of the map is multiplied by -1 since we increase the number of edges without changing the number of white vertices. For $i\geq1$, the operator $x_i^2\frac{\partial}{\partial x_i}$ allows to choose a $\circ^1$-corner in a face of color $i$ and to increase the degree $\deg_{\circ^1}$ of this face by 1. Finally, we have two ways to distinguish a side of the added edge and this compensated by the factor $1/2^{|\lambda|-\cc(M)}$.

\cref{eq Y+} can be obtained in a similar way.
\end{proof}
\end{lem}


We deduce the following proposition.

\begin{prop}\label{prop F}
Let $\lambda=2^r1^s$ be a  2-column partition. One has
$$\Frho(\bfx)=\sum_{0\leq j\leq r}(r)_jX_{+}^{r+s-j}Y_{+}^{r-j} \Bj(\bfxy,-\alpha,-\alpha)_{\Big|\bfy=\bfx}.$$

\begin{proof}
We start by noticing that substituting $\bfx$ by $\bfxy$ in \cref{eq B}, specializing $u_1=u_2=-\alpha$ and developing the power-sum functions as in \cref{eq p x}, we get
$$\Bj(\bfxy,-\alpha,-\alpha)=\frac{1}{j!}\sum_{(M_\circ,\mathcal{C})}(-1)^{|\ver(M_\circ)|}\alpha^{|\ver(M_\circ)|-\cc(M_\circ)}b^{\wvt(M_\circ)}\prod_f \left(x_{\mathcal{C}(f)}y_{\mathcal{C}(f)}\right)^{\deg(f)},$$
where the sum runs over colored labelled rooted bipartite maps $(M_\circ,\mathcal{C})$ of size $j$ which are bipartite in the colors $\circ^1/\circ^2$. Realizing the three first steps detailed in \cref{ssec decorated labelled maps} on each map $M_\circ$, the previous equation can be rewritten as follows
\begin{align}\label{eq B prop F}
\binom{r}{j}\Bj(\bfxy,-\alpha,-\alpha)
&=\frac{1}{j!}\sum_{\text{colored rooted}\atop{\text{maps }(M,\mathcal{C})}}(-1)^{|\wv(M)|}\alpha^{|\wv(M)|-\cc(M)}b^{\wvt(M)}\kappa(M,\mathcal{C)}\\
&=\frac{1}{j!}\sum_{\text{colored side-}\atop{\text{marked maps }(M,\mathcal{C})}}\frac{(-1)^{|\wv(M)|}}{2^{2j-\cc(M)}}\alpha^{|\wv(M)|-\cc(M)}b^{\wvt(M)}\kappa(M,\mathcal{C)}\nonumber,
\end{align}
where the two sums run over $\lambda$-injectively decorated maps $(M,\mathcal{C})$ of size $2j$, which are bipartite in the two colors white and black, and such that all black vertices have degree 2. 
In order to add black-leaf edges as explained in the last step at the end of \cref{ssec decorated labelled maps}, we apply the operators $X_+$ and $Y_+$. Fix a map $(M,\mathcal{C})$ as in \cref{eq B prop F}. Using \cref{lem add black leaves}, we get that  
$$X_{+}^{r+s-j}Y_{+}^{r-j}\frac{1}{2^{2j-\cc(M)}}\kappa(M,\mathcal{C})=\sum_{(N,\mathcal{D})}\frac{1}{2^{|\lambda|-\cc(N)}}\kappa(N,\mathcal{D}),$$
where the sum runs over colored
$\lambda$-injectively decorated side-marked bipartite maps $(N,\mathcal{D})$ of size $|\lambda|$, such that $M$ is obtained from $N$ by deleting all black-leaf edges and the coloring $\mathcal D$ is inherited from $\mathcal{C}$. This finishes the proof of the proposition.
\end{proof}
\end{prop}

\subsection{End of the proof of the main result}
Since the generating function $\Frho$ is obtained from the functions $\Bj(\bfxy,-\alpha,-\alpha)$ by applying the operators $X_+$ and $Y_+$ (see \cref{prop F}) and since this function has an expression in terms of 1-column monomial functions (see \cref{thm bipartite maps}), we should understand the action of the operators $X_+$ and $Y_+$ on 1-column monomial functions in $\bfxy$. This is given by the following lemma.

\begin{lem}\label{lem X+ Y+}
Let $r,s$ and $j$ be three non-negative integers satisfying $j\leq r$. Then 
$$X_{+}^{r+s-j}Y_{+}^{r-j}m_{1^j}(\mathbf{xy})_{\Big|\bfy=\bfx}=\sum_{j\leq i \leq r}\binom{i}{j}\binom{2(r-i)+s}{r-i}(r-j)!(r+s-j)!m_{2^i,1^{2(r-i)+s}}.$$
\begin{proof}
For any subset  $\beta\subset\mathbb N^*$ and an alphabet of variables $\mathbf{x}=(x_1,x_2,...)$ we define 
$$x^{\beta}:=\prod_{i\in\beta}x_i,\hspace{0.5cm}\text{and} \hspace{0.5cm}x^{2\beta}:=\prod_{i\in\beta}x_i^2.$$
With this notation, we write
\begin{equation}\label{eq column monomials}
  m_{1^j}(\mathbf{xy})=\sum_{|\beta|=j}x^{\beta}y^{\beta}.  
\end{equation}

We start by noticing that for every subset $\beta\subset \mathbb{N^*}$, one has
$$X_+x^\beta=\sum_{i\notin \beta}x^{\beta\cup \{i\}}.$$
Hence, for every $\beta$ of size $j$ we have
\begin{equation}\label{eq X+Y+}
X_+^{r+s-j}Y_+^{r-j}x^\beta y^\beta=(r+s-j)!(r-j)!\sum_{\gamma,\delta_1,\delta_2}x^{\gamma\cup\delta_1}y^{\gamma\cup\delta_2}    
\end{equation}

where the sum is taken over $\gamma,\delta_1,\delta_2\subset \mathbb N^*$, satisfying
\begin{itemize}
    \item $|\gamma|+|\delta_1|=r+s$ and $|\gamma|+|\delta_2|=r$.
    \item $\beta\subset\gamma$.
    \item The three sets $\delta_1$, $\delta_2$ and $\gamma$ are disjoint.
\end{itemize}  
Taking the sum over sets $\beta$ of size $j$ and setting $\bfy=\bfx$ gives
\begin{align*}
    X_+^{r+s-j}Y_+^{r-j} m_{1^j}(\mathbf{xy})_{\Big|\bfy=\bfx}
    &=(r+s-j)!(r-j)!\sum_{|\beta|=j}\sum_{\gamma,\delta_1,\delta_2}x^{2\gamma}x^{\delta_1\cup\delta_2}\\
    &=(r-j)!(r+s-j)!\sum_{j\leq i\leq r}\sum_{\gamma,\delta}C_j(\gamma,\delta) x^{2\gamma}x^\delta ,\\
\end{align*}
where the second sum in the last line is taken over disjoint sets $\gamma$ and $\delta$ such that $|\gamma|=i$ and $|\delta|=2(r-i)+s$, and where $C_j(\gamma,\delta)$ is the number of triplets of sets $(\beta,\delta_1,\delta_2)$ of respective sizes  $j,r+s-i$ and $r-i$ and such that $\beta\subset\gamma$ and $(\delta_1,\delta_2)$ is a partition of $\delta$. Hence

\begin{align*}
C_j(\gamma,\delta)
=\binom{2(r-i)+s}{r-i}\binom{i}{j},
\end{align*}
and this finishes the proof of the lemma.
\end{proof}
\end{lem}

We now deduce the proof of the main theorem.

\begin{proof}[Proof of the main theorem]
Fix a MON $\rho$. Using \cref{thm bipartite maps} and \cref{prop F} we get that 
$$\Frho=\sum_{0\leq j\leq r} (r)_j(\alpha+j-1)_jX_{+}^{r+s-j}Y_{+}^{r-j}m_{1^j}(\mathbf{xy})_{\Big|\bfy=\bfx}.$$
On the other hand, using \cref{lem X+ Y+} we get that
\begin{align*}
    \sum_{0\leq j\leq r}& (r)_j(\alpha+j-1)_jX_{+}^{r+s-j}Y_{+}^{r-j}m_{1^j}(\mathbf{xy})_{\Big|y_i=x_i} \\    
    &=\sum_{0\leq j\leq r}(r)_j(\alpha+j-1)_j\sum_{j\leq i\leq r}\binom{i}{j}\binom{2(r-i)+s}{r-i}(r-j)!(r+s-j)!m_{2^i,1^{2(r-i)+s}}(\bfx)\\
    &=\sum_{0\leq i\leq r}r!\binom{2(r-i)+s}{r-i}m_{2^i,1^{2(r-i)+s}}(\bfx)\sum_{0\leq j\leq i}\binom{i}{j}(\alpha+j-1)_j(r+s-j)!\\
    &=\sum_{0\leq i\leq r}(r)_i (2(r-i)+s)!m_{2^i,1^{2(r-i)+s}}(\bfx)\sum_{0\leq j\leq i}\frac{\binom{i}{j}(\alpha+j-1)_j(r+s-j)!}{(r+s-i)!}.
\end{align*}
The second sum in the last line can be rewritten as follows
\begin{align*}
  (-1)^ii!\sum_{0\leq j\leq i}\binom{-\alpha}{j}\binom{-(r+s-i+1)}{i-j}
  &=(-1)^ii!\binom{-(\alpha+r+s-i+1)}{i}\\
  &=(\alpha+r+s)_i.
\end{align*}

\noindent Hence, we obtain that 
\begin{align*}\label{Eq 2col}
\Frho=\sum_{i=1}^r(r)_i(\alpha+r+s)_i(2(r-i)+s)!m_{2^i1^{2(r-i)+s}}(\mathbf{x}).
\end{align*}
Using \cref{eq 2-column Jack}, we get that 
$\Frho=J_\lambda^{(\alpha)}.$
\end{proof}

\section{Application: \texorpdfstring{\cref{conj injective}}{} for the low degree coefficients in \texorpdfstring{$\alpha$}{}}\label{sec applications}

We recall that $\alpha$ is the Jack parameter, and $b$ is the parameter related to $\alpha$ by $\alpha=b+1$. Let $R$ be a field.
We denote by $R (\alpha)$ the field of rational functions in $\alpha$ with coefficients in $R$.
For  $f\in R (\alpha)$ and an integer $m$, we write $f=O(\alpha^m)$ if the rational function $\alpha^{-m}\cdot f$ has no pole in 0.
The main result of this section establishes the cases of $\alpha=0$ and the coefficient $[\alpha]$ in \cref{conj injective}.
\begin{thm}\label{thm conj 0 1}
Let $\lambda$ be a partition of size $n$.  Then, there exists a \textup{MON} $\rho$ such that
$$J_\lambda^{(\alpha)}(\mathbf{p})=\sum_{M\in\ID_n(\lambda)}(-1)^{n-|\wv(M)|}\alpha^{|\wv(M)|-\cc(M)}b^{\vartheta^\rho(M)}p_{\type(M)}+O(\alpha^2).$$

\end{thm}

Let $\lambda^1,...,\lambda^r$ be a family of partitions. We denote by $\bigoplus_{1\leq i\leq r}\lambda^i$ its \textit{entry-wise sum} defined by $\left(\bigoplus_{1\leq i\leq r}\lambda^i\right)_j=\sum_{1\leq i\leq r}\lambda^i_j$, for every $j\geq 1$. In particular, a partition $\lambda$ can be written as the entry-wise sum of its columns;
$\lambda=\oplus_{1\leq i\leq \lambda_1}1^{\lambda^t_i},$
where $\lambda^t$ denotes the conjugate partition of $\lambda$.

The following theorem is a particular case of the strong factorization property of Jack polynomials due to Do\l{}e\k{}ga and Féray \cite[Theorem 1.4]{DolegaFeray2017}.

\begin{thm}\cite{DolegaFeray2017}\label{thm cumulants}
Let $\lambda^1$, $\lambda^2$ and $\lambda^3$ be three partitions.
Then 
$$J^{(\alpha)}_{\lambda^1\oplus\lambda^2\oplus\lambda^3}-J^{(\alpha)}_{\lambda^1\oplus\lambda^2}J^{(\alpha)}_{\lambda^3}-J^{(\alpha)}_{\lambda^1\oplus\lambda^3}J^{(\alpha)}_{\lambda^2}-J^{(\alpha)}_{\lambda^2\oplus\lambda^3}J^{(\alpha)}_{\lambda^1}+2J^{(\alpha)}_{\lambda^1}J^{(\alpha)}_{\lambda^2}J^{(\alpha)}_{\lambda^3}=O(\alpha^2).$$
\end{thm}
For $\alpha=0$, we have the following expression for Jack polynomials , see \cite[Proposition 7.6]{Stanley1989}.

\begin{equation}\label{eq alpha=0}
J_\lambda^{(0)}=\prod _{1\leq i\leq\lambda_1}J_{\lambda^t_i}^{(0)}.
\end{equation}
Using the strong factorization property we can generalize this result as follows; the coefficient of $\alpha^{r-1}$ in the expansion of a Jack polynomial as a polynomial in $\alpha$, can be obtained using only Jack polynomials indexed by partitions with less than $r$ columns. We prove here this result for $r=2$ using \cref{thm cumulants}.

\begin{cor}
Let $\lambda$ be partition with 3 columns or more. Then,

\begin{equation}\label{eq deg1 alpha}
[\alpha]J_\lambda^{(\alpha)}=\sum_{1\leq i<j\leq \lambda_1}\left([\alpha]J_{1^{\lambda^t_i}\bigoplus1^{\lambda^t_j}}^{(\alpha)}\right)\prod_{k\neq,i,j}J_{1^{\lambda^t_k}}^{(0)}.    
\end{equation}

\begin{proof}
We prove the result by induction on the number of columns of $\lambda$. If $\lambda$ has 3 columns then the result is direct consequence of \cref{thm cumulants} and the fact that Jack polynomials indexed by 1-column partitions are independent from $\alpha$ (see \cref{eq column Jack}).
Let $\lambda$ be a partition containing more than 3 columns. We denote by $1^{\lambda^t_1}$ and $1^{\lambda^t_2}$ the two first columns of $\lambda$ and by $\mu$ the partition obtained by taking the entry-wise sum of the other columns. Hence we have $\lambda=1^{\lambda^t_1}\oplus1^{\lambda^t_2}\oplus\mu$.
From \cref{thm cumulants} we get that 
\begin{align*}
[\alpha]J^{(\alpha)}_\lambda
&=[\alpha]\left[J^{(\alpha)}_{1^{\lambda^t_1}\oplus1^{\lambda^t_2}}J^{(\alpha)}_{\mu}+J^{(\alpha)}_{1^{\lambda^t_1}\oplus\mu}J^{(0)}_{1^{\lambda^t_1}}+J^{(\alpha)}_{1^{\lambda^t_2}\oplus\mu}J^{(0)}_{1^{\lambda^t_1}}-2J^{(0)}_{1^{\lambda^t_1}}J^{(0)}_{1^{\lambda^t_2}}J^{(\alpha)}_{\mu}\right]\\
&=\left([\alpha]J^{(\alpha)}_{1^{\lambda^t_1}\oplus1^{\lambda^t_2}}\right)J^{(0)}_{\mu}+J^{(0)}_{1^{\lambda^t_1}\oplus1^{\lambda^t_2}}\left([\alpha]J^{(\alpha)}_{\mu}\right)+\left([\alpha]J^{(\alpha)}_{1^{\lambda^t_1}\oplus\mu}\right)J^{(0)}_{1^{\lambda^t_1}}\\
&\hspace{3.5cm}+\left([\alpha]J^{(\alpha)}_{1^{\lambda^t_2}\oplus\mu}\right)J^{(0)}_{1^{\lambda^t_1}}-2J^{(0)}_{1^{\lambda^t_1}}J^{(0)}_{1^{\lambda^t_2}}\left([\alpha]J^{(\alpha)}_{\mu}\right).    
\end{align*}
For each $1\leq i< j\leq \lambda_1$, we prove that the term $$\left([\alpha]J_{1^{\lambda^t_i}\bigoplus1^{\lambda^t_j}}^{(\alpha)}\right)\prod_{k\neq,i,j}J_{1^{\lambda^t_k}}^{(0)}$$
appears with coefficient 1 in the last sum by applying the induction hypothesis on the partitions $1^{\lambda^t_1}\oplus \mu$, $1^{\lambda^t_2}\oplus \mu$ and $\mu$ and then using  \cref{eq alpha=0}; we distinguish the three cases 
\begin{itemize}
    \item $(i,j)=(1,2)$,
    \item $i\leq 2$ and $j>2$,
    \item $2<i<j$.\qedhere
\end{itemize}

\end{proof}
\end{cor}

The following lemma will be useful in the proof of \cref{prop F}.
\begin{lem}\label{lem star maps}
Fix a partition $\lambda$ and a $\lambda$-injectively decorated map $M$, such that $|\mathcal{V}_\circ(M)|-\cc(M)=0$. Then, $M$ is orientable and then $\vartheta(M)=0$ for every statistic of non-orientability $\vartheta$. 
\begin{proof}
Since the number of white vertices of $M$ is equal to the number of connected component, then $M$ has exactly one white vertex in each connected component. On the other hand, since $M$ is $\lambda$-injectively decorated then the underlying graph does not have multiple edges. Hence, all black vertices of $M$ has degree 1, and $M$ is orientable.
\end{proof}
\end{lem}

It can be shown that the right-hand side of \cref{conj injective} also satisfies the strong factorization of Féray and Do\l{}e\k{}ga. We prove here the two first equations of this property which will be useful for the proof of \cref{thm conj 0 1}.
\begin{prop}
For every non-orientability statistic $\vartheta$, we have 
\begin{equation}\label{eq F alpha=0}
  F^{(0)}_{\lambda,\vartheta}=\prod _{1\leq i\leq\lambda_1}F_{\lambda^t_i,\vartheta}^{(0)},  
\end{equation}

and 
\begin{equation}\label{eq F alpha}
[\alpha]\F=\sum_{1\leq i<j\leq \lambda_1}\left([\alpha]F_{1^{\lambda^t_i}\bigoplus1^{\lambda^t_j},\vartheta}^{(\alpha)}\right)\prod_{k\neq,i,j}F_{1^{\lambda^t_k},\vartheta}^{(0)}.    
\end{equation}
\begin{proof}

When $\alpha=0$, only the maps $M$ such that $|\mathcal{V}_\circ(M)|-\cc(M)=0$ appear in the sum defining $\F$, \textit{i.e} maps having exactly one white vertex in each connected component. Hence, in such a map the edges of one connected component are labelled by boxes in the same column of $\lambda$. We deduce that for such a map $M$ (necessarily orientable by \cref{lem star maps}), we have 
\begin{equation}\label{eq prop F}
(-1)^{n-|\wv(M)|}p_{\type(M)}=\prod_{1\leq i\leq \lambda_1}(-1)^{|M_i|-|\wv(M_i)|}p_{\type(M_i)},
\end{equation}
where for each $i$, $M_i$ denotes the collection of connected components of $M$, whose edges are decorated by boxes in the column $i$ of $\lambda.$
This finishes the proof of \cref{eq F alpha=0}.

When we consider the coefficient  $[\alpha]\F$, only maps $M$ such that~$|\wvt(M)|-\cc(M)=1$ contribute to the sum (this is a consequence of \cref{lem star maps}).  Each connected component of such a map contains exactly one white vertex, except for one that contains two white vertices. It is easy to check that the edges of such a connected component are labelled by boxes in two different columns. Using multiplicativity arguments as in \cref{eq prop F}, we deduce \cref{eq F alpha}.
\end{proof}

\end{prop}

We deduce the proof of \cref{main thm}.
\begin{proof}[Proof of \cref{thm conj 0 1}]
We use the fact that the coefficients $[\alpha^0]$ and $[\alpha]$ of Jack polynomials and of the generating series $\F$ satisfy the same equations (\cref{eq alpha=0,eq deg1 alpha} on one hand and \cref{eq F alpha=0,eq F alpha} on the other hand) and that in these equations only functions $J^{(\alpha)}_\lambda$ and $\F$ indexed by 1- and 2- column partitions are involved. But we know from \cref{main thm} that there exist a MON $\rho$ such that \cref{conj injective} holds in these cases.
\end{proof}

\appendix 
\section{Proof of \texorpdfstring{\cref{thm CD20}}{}}\label{appendix}

Following \cite{ChapuyDolega2022}, we introduce the differential operators
$$A_1:=p_1/\alpha,\qquad A_2:=\sum_{i\geq1}p_{i+1}\frac{i\partial}{\partial p_{i}},$$ 
$$A_3:=(1+b)\sum_{i,j\geq1}p_{i+j+1}\frac{ij\partial^2}{\partial p_{i}\partial p_{j}}+\sum_{i,j\geq1}p_i p_j\frac{(i+j)\partial}{\partial p_{i+j-1}}+b\sum_{i\geq1}p_{i+1}\frac{i\partial}{\partial p_i},$$
on the algebra $\mathbb{Q}(b)[p_1,p_2,..]$.
In the following proposition, we give a combinatorial interpretation of these operators.

\begin{prop}\label{prop add edges}
Fix a $\MON$ $\rho$ and a labelled side-marked bipartite map $M$ of size $n$. 
Then,
$$A_i \left[\frac{b^{\wvt(M)}}{2^{|M|-\cc(M)}\alpha^{\cc(M)}}p_{\type(M)}\right]
\\=\sum_{\widetilde{M}}\frac{b^{\wvt(\widetilde{M})}}{2^{|\widetilde{M}|-\cc(\widetilde{M})}\alpha^{\cc(\widetilde{M})}}p_{\type(\widetilde{M})}, \quad \text{for } 1\leq i\leq3,    $$

where the sum is taken over labelled side-marked maps $\widetilde{M}=M\cup \{e\}$ obtained by adding a side-marked edge $e$ of label $n+1$ to the map $M$, such that:
\begin{enumerate}[label={\textup{(}\arabic*\textup{)}}]
    \item if $i=1$, then $e$ is a disconnected edge.
    \item if $i=2$, then $e$ is a leaf-edge connecting a white corner of $M$ to an isolated black vertex (or equivalently a black corner of $M$ to an isolated white vertex).
    \item if $i=3$, then $e$ connects two corners  of different colors in $M$.
\end{enumerate}

\begin{proof}
We start by observing that the action of the operator $\frac{i\partial}{\partial p_i}$ on $p_{\type(M)}$ can be interpreted as choosing a white corner (or equivalently a black corner) in a face of degree $i$ in $M$ (we recall that such a face contain $i$ corners of each color).

In the case of item (\textit{1}), we add a face of size 1 and the number of connected components increases by 1, hence the multiplication by $p_1/\alpha$.

In the case of item (\textit{2}), we choose a white corner in a face $f$ and we increase the degree of $f$ by 1.
In these two cases, the $b$-weight of the map does not change by adding the edge $e$.

We now focus on item (\textit{3}). We distinguish three cases.
\begin{itemize}
    \item We add the edge $e$ between two corners $c_1$ and $c_2$ which are incident to two different faces $f_1$ and $f_2$ of respective sizes $i$ and $j$, to form a face of size $i+j+1$. Let us show that this case corresponds to the first term of the operator $A_3$. Let $\tilde{e}$ denote the edge obtained by twisting $e$, see \cref{ssec stat}. 
    
    If the two faces lie in the same connected component of $M$, then $e$ is a handle and  by definition of a MON  $$\rho(M\cup\{e\},e)+\rho(M\cup\{\tilde e\},\tilde e)=1+b,$$
    and this explains the factor $1+b$ in the first terms of $A_3$.
    
    If the two faces lie in two different connected components, then $e$ is bridge and $$\rho(M\cup\{e\},e)=\rho(M\cup\{\tilde e\},\tilde e).$$ In this case, the factor $1+b=\alpha$ in the first term of $A_3$ is related to the fact that the number of connected components of $M$ decreases by 1.
    
    \item The added edge $e$ is a diagonal between two corners incident to a face of degree $i+j-1$, which splits the face into two faces of respective degrees $i$ and $j$. Then $\rho(M\cup\{e\},e)=1$ and both the $b$- and the $\alpha$-weight of the map are unchanged.
    
    \item The edge $e$ a twist added on a face of degree $i\geq1$ to obtain a face of degree $i+1$. Then  $\rho(M\cup\{e\},e)=b$.
\end{itemize}

Note that in the cases of items (\textit{2}) and (\textit{3}), we have each time two ways to distinguish a side on the added edge, so that the map obtained is side-marked. This explains the factor $2^{|M|-\cc(M)}$ which appears in the denominator. 
\end{proof}
\end{prop}

We deduce the following theorem.
\begin{thm}\label{prop B}
Fix a $\MON$ $\rho$. For every $n\geq 0$, we have
\begin{equation}\label{eq decomposition}
    (n+1)B^{(\alpha)}_{n+1,\rho}(\bfx,u_1,u_2)=(A_3+(u_1+u_2)A_2+u_1u_2A_1)B^{(\alpha)}_{n,\rho}(\bfx,u_1,u_2).
\end{equation}
\begin{proof}
We recall that by definition of the function $\B$ has the following expression 
\begin{equation*}
\B(\bfx,u_1,u_2)= \sum_{M}\frac{1}{n!2^{n-\cc(M)}}u_1^{|\wv(M)|}u_2^{|\bv(M)|}\alpha^{-\cc(M)}b^{\wvt(M)}p_{\type(M)}(\bfx),    
\end{equation*}
where the sum is taken over labelled side-marked maps of size $n$.
The theorem is then a direct consequence of \cref{prop add edges}.
\end{proof}
\end{thm}
\begin{rmq}
The previous theorem is a variant of the decomposition equation of $k$-constellation established in \cite[Theorem 3.10]{ChapuyDolega2022}; the case of bipartite maps that we consider here corresponds to $k=2$ and to the specialization $\mathbf{q}:=(1,0,0,..)$ with the notation of \cite{ChapuyDolega2022} (see also the proof of \cref{thm CD20} below).
\end{rmq}
We now deduce the proof of \cref{thm CD20}.

\begin{proof}[Proof of \cref{thm CD20}]
Let $\mathbf{q}=(q_1,q_2,...)$ be an additional sequence of variables.
We consider the following function introduced in \cite{ChapuyDolega2022}.
$$\tau_b^{(2)}(t,\mathbf{p},\mathbf{q},u_1,u_2)=\sum_{j\geq0}t^n\sum_{\xi\vdash j}\frac{J^{(\alpha)}_\xi(\mathbf{p})J^{(\alpha)}_\xi(\mathbf{q})J^{(\alpha)}_\xi(\underline{u_1)}J^{(\alpha)}_\xi(\underline{u_2)}}{j^{(\alpha)}_\xi}.$$
From \cref{prop B}, we obtain that $\sum_{n\geq 0} t^n\B(\bfx,u_1,u_2)$ satisfies the differential equation \cite[Eq.(16)]{ChapuyDolega2022} for $k=2$, $m=1$ and $\mathbf{q}=\delta_1:=(1,0,0,...)$. Since this equation fully characterizes the functions $\B$, then
using \cite[Theorem 5.10]{ChapuyDolega2022}, we get that 
$$\sum_{n\geq 0} t^n\B(\bfx,u_1,u_2)=\tau_b^{(2)}(t,\mathbf{p},\delta_1,u_1,u_2).$$
We conclude using the fact that $J^{(\alpha)}_\lambda(\delta_1)=1$ (see \cref{Jack formula}).
\end{proof}

\noindent {\bf Acknowledgements.} The author is very grateful to his advisors Valentin Féray and Guillaume Chapuy for several interesting discussions about Jack polynomials and maps enumeration.

\bibliographystyle{alpha}
\bibliography{biblio2}

\end{document}